\documentclass[10pt]{amsart}
\usepackage{latexsym}   
\usepackage{amssymb}    
\usepackage{amsmath}    
\usepackage{amsthm}     
\usepackage[all]{xy}
\usepackage{hyperref}
\newcommand{\pa}{\partial}

\newcommand{\Ga}{\Gamma}\newcommand{\del}{\delta}

\newcommand{\la}{\lambda}\newcommand{\La}{\Lambda}\newcommand{\om}{\omega}
\newcommand{\Om}{\Omega}
\newcommand{\ti}{\tilde}

\renewcommand{\thefootnote}

\theoremstyle{definition}

\theoremstyle{remark}

\numberwithin{equation}{section}

\title[ Multiply conjugate systems containing deformations of quadrics]
{ Multiply conjugate systems containing deformations of quadrics}

\author[  Ion I. Dinc\u{a}]{Ion I. Dinc\u{a}}
\address{Faculty of Mathematics and Informatics,
University of Bucharest,  14 Academiei Str., 010014, Bucharest,
Romania}
 \email{dinca@gta.math.unibuc.ro}
\thanks{Supported by the University of Bucharest}

\begin{document}


\begin{abstract}
We provide a generalization of Bianchi's triply conjugate systems
containing a family of deformations of $2$-dimensional quadrics
together with its B\"{a}cklund transformation to higher
dimensions.
\end{abstract}

\maketitle

\tableofcontents \pagenumbering{arabic}

\section{Introduction}

At the end of the XIX$^{\mathrm{th}}$ century Weingarten produced
the first examples of triply orthogonal systems (Lam\'{e}
families) containing a family of {\it constant Gau\ss\ curvature}
(CGC) $K=-1$ surfaces. Bianchi extended his work by allowing the
(negative) CGC to vary within its family of surfaces and by
developing the {\it B\"{a}cklund} (B) transformation for such
triply orthogonal systems. With the development of the theory of
deformations of quadrics at the beginning of the
XX$^{\mathrm{th}}$ century Bianchi introduced triply conjugate
systems containing a family of deformations of quadrics together
with its B transformation as the natural metric-projective
generalization of Weingarten's triply orthogonal systems (see
Darboux \cite{D2} for triply conjugate systems, Green \cite{G} for
projective differential geometry of triple systems of surfaces and
Bianchi (\cite{B4}, Vol {\bf 4},(143)\&(146)) for triply conjugate
systems containing deformations of quadrics).

According to Bianchi such triply conjugate systems must induce on
the deformations of quadrics the conjugate system common to the
quadric and its deformation, thus generalizing lines of curvature
on CGC surfaces (surfaces of triply orthogonal systems cut each
other along lines of curvature).

Note that Peterson's $1$-dimensional family of deformations of
quadrics (see \cite{P}) a-priori does not provide a triply
conjugate system, although we have the conjugate system common to
the quadric and its deformation property for each deformation.

In what concerns deformations of higher dimensional quadrics the
first results are those of Cartan \cite{C} concerning the
deformation problem of space forms in space forms; in particular
he introduced the minimal co-dimension $(n-1)$ for the deformation
problem of $n$-dimensional space forms in space forms, the
exteriorly orthogonal forms tool (naturally appearing from the
Gau\ss\ equations), their canonical form (given by lines of
curvature on the deformation) and the $n(n-1)$ functions of one
variable dimensionality of the space
of such deformations. These have been extended on one hand (upon a suggestion
from S. S. Chern and using a result due to Moore on the Chebyshev
coordinates on deformations of $\mathbb{H}^n(\mathbb{R})$ in
$\mathbb{R}^{2n-1}$; these are lines of curvature and thus in
bijective correspondence with such deformations) to the B
transformation of $\mathbb{H}^n$ in $\mathbb{R}^{2n-1}$ in
Tenenblat-Terng \cite{TT} (Terng also developed the {\it Bianchi
Permutability Theorem} (BPT) for this B transformation in
\cite{T}) and on the other hand to deformations in
$\mathbb{R}^{2n-1}$ of quadrics in $\mathbb{R}^{n+1}$ or in
$\mathbb{R}^n\times(i\mathbb{R})$ with positive definite linear
element in Berger, Bryant and Griffiths \cite{BBG} (again we have
the minimal co-dimension $(n-1)$ and the $n(n-1)$ functions of one
variable dimensionality of the space
of such deformations).

Note also that multiply orthogonal systems are present in the
current literature of integrable systems (see for example
Terng-Uhlenbeck \cite{TU} and its references).

Thus the natural question of completing deformations of higher
dimensional quadrics and their B transformation to multiply
conjugate systems containing deformations of higher dimensional
quadrics and their B transformation arises. Just as the theory of
deformation of  $2$-dimensional quadrics admits discretization via
the iteration of the B transformation (moving M\"{o}bius
configurations in Bianchi's denomination; see Bobenko-Pinkall
\cite{BP} for discrete deformations of the $2$-dimensional
pseudo-sphere), a similar approach should give discrete multiply
conjugate systems containing deformations of higher dimensional
quadrics (note that Weingarten discovered his triply orthogonal
systems of CGC $-1$ surfaces by the iteration of an infinitesimal
transformation of CGC $-1$ surfaces, so he essentially used the
discrete version to find the differential version).

According to the principles laid down by Bianchi, if we restrict
$(n-1)$ parameters in the multiply conjugate systems to constants
so as to obtain a deformation of a quadric, then the remaining $n$
parameters should provide the conjugate system common to the
quadric and its deformation (a-priori the quadric may vary with
the $(n-1)$ constants, but we restrict our discussion only to
deformations of a fixed quadric).

Thus to obtain multiply conjugate systems containing deformations
of quadrics one must extend the differential system for
deformations of $n$-dimensional quadrics in $\mathbb{C}^{2n-1}$
via isothermal-conjugate system on the considered quadric and the
conjugate system common to the quadric and its deformation
according to Bianchi's principles and Cartan's exterior
differential calculus.

For more details on the (classical) theory of deformations of
(higher dimensional) quadrics we refer the reader to one of our
previous notes concerning Bianchi's B\"{a}cklund transformation
for higher dimensional quadrics.

All computations are local and assumed to be valid on their open
domain of validity without further details; all functions have the
assumed order of differentiability and are assumed to be
invertible, non-zero, etc when required (for all practical
purposes we can assume all functions to be analytic).

\section{Confocal quadrics in canonical form}

Consider the complexified Euclidean space
$$(\mathbb{C}^m,<.,.>),\
<x,y>:=x^Ty,\ |x|^2:=x^Tx,\ x,y\in\mathbb{C}^m$$ with standard
basis $\{e_j\}_{j=1,...,m},\ e_j^Te_k=\del_{jk}$.

Isotropic (null) vectors are those vectors $v$ of length $0\
(|v|^2=0)$; since most vectors are not isotropic we shall call a
vector simply vector and we shall only emphasize isotropic when
the vector is assumed to be isotropic. The same denomination will
apply in other settings: for example we call quadric a
non-degenerate quadric (a quadric projectively equivalent to the
complex unit sphere).

A quadric $x\subset\mathbb{C}^{n+1}$ is given by the quadratic
equation
$Q(x):=\begin{bmatrix}x\\1\end{bmatrix}^T\begin{bmatrix}A&B\\B^T&C\end{bmatrix}
\begin{bmatrix}x\\1\end{bmatrix}=x^T(Ax+2B)+C=0,\ A=A^T\in\mathbf{M}_{n+1}(\mathbb{C}),\
B\in\mathbb{C}^{n+1},\ C\in\mathbb{C},\
\begin{vmatrix}A&B\\B^T&C\end{vmatrix}\neq 0$.

A metric classification of all (totally real) quadrics in
$\mathbb{C}^{n+1}$ requires the notion of {\it symmetric Jordan}
(SJ) canonical form of a symmetric complex matrix. The symmetric
Jordan blocks are: $J_1:=0=0_{1,1}\in\mathbf{M}_1(\mathbb{C}),\
J_2:=f_1f_1^T\in\mathbf{M}_2(\mathbb{C}),\
J_3:=f_1e_3^T+e_3f_1^T\in\mathbf{M}_3(\mathbb{C}),\ J_4:=f_1\bar
f_2^T+f_2f_2^T+\bar f_2f_1^T\in\mathbf{M}_4(\mathbb{C}),\ J_5:=
f_1\bar f_2^T+f_2e_5^T+e_5f_2^T+\bar
f_2f_1^T\in\mathbf{M}_5(\mathbb{C}),\ J_6:= f_1\bar f_2^T+f_2\bar
f_3^T+f_3f_3^T+\bar f_3f_2^T+\bar
f_2f_1^T\in\mathbf{M}_6(\mathbb{C})$, etc, where
$f_j:=\frac{e_{2j-1}+ie_{2j}}{\sqrt{2}}$ are the standard
isotropic vectors (at least the blocks $J_2,\ J_3$ were known to
the classical geometers). Any symmetric complex matrix can be
brought via conjugation with a complex rotation to the symmetric
Jordan canonical form, that is a matrix block decomposition with
blocks of the form $a_jI_p+J_p$; totally real quadrics are
obtained for eigenvalues $a_j$ of the quadratic part $A$ defining
the quadric being real or coming in complex conjugate pairs $a_j,\
\bar a_j$ with subjacent symmetric Jordan blocks of same dimension
$p$. Just as the usual Jordan block $\sum_{j=1}^pe_je_{j+1}^T$ is
nilpotent with $e_{p+1}$ cyclic vector of order $p$, $J_p$ is
nilpotent with $\bar f_1$ cyclic vector of order $p$, so we can
take square roots of SJ matrices without isotropic kernels
($\sqrt{aI_p+J_p}:=\sqrt{a}\sum_{j=0}^{p-1}(^{\frac{1}{2}}_j)a^{-j}J_p^j,\
a\in\mathbb{C}^*,\ \sqrt{a}:=\sqrt{r}e^{i\theta}$ for
$a=re^{2i\theta},\ 0<r,\ -\pi\le 2\theta<\pi$), two matrices with
same SJ decomposition type (that is $J_p$ is replaced with a
polynomial in $J_p$) commute, etc.

The confocal family $\{x_z\}_{z\in\mathbb{C}}$ of a quadric
$x_0\subset\mathbb{C}^{n+1}$ in canonical form (depending on as
few constants as possible) is given in the projective space
$\mathbb{C}\mathbb{P}^{n+1}$ by the equation
$Q_z(x_z):=\begin{bmatrix}x_z\\1\end{bmatrix}^T(\begin{bmatrix}A&B\\B^T&C\end{bmatrix}^{-1}-z
\begin{bmatrix}I_{n+1}&0\\0^T&0\end{bmatrix})^{-1}\begin{bmatrix}x_z\\1\end{bmatrix}=0$,
where

$\bullet\ \ A=A^T\in\mathbf{GL}_{n+1}(\mathbb{C})$ SJ,
$B=0\in\mathbb{C}^{n+1},\ C=-1$ for {\it quadrics with center}
(QC),

$\bullet\ \ A=A^T\in\mathbf{M}_{n+1}(\mathbb{C})$ SJ,
$\ker(A)=\mathbb{C}e_{n+1},\ B=-e_{n+1},\ C=0$ for {\it quadrics
without center} (QWC) and

$\bullet\ \ A=A^T\in\mathbf{M}_{n+1}(\mathbb{C})$ SJ,
$\ker(A)=\mathbb{C}f_1,\ B=-\bar f_1,\ C=0$ for {\it isotropic
quadrics without center} (IQWC).

From the definition one can see that the family of quadrics
confocal to $x_0$ is the adjugate of the pencil generated by the
adjugate of $x_0$ and Cayley's absolute
$C(\infty)\subset\mathbb{C}\mathbb{P}^n$ in the hyperplane at
infinity; since Cayley's absolute encodes the Euclidean structure
of $\mathbb{C}^{n+1}$ (it is the set invariant under rigid motions
and homotheties of
$\mathbb{C}^{n+1}:=\mathbb{CP}^{n+1}\backslash\mathbb{CP}^n$) the
mixed metric-projective character of the confocal family becomes
clear.

For QC $\mathrm{spec}(A)$ is unambiguous (does not change under
rigid motions) but for (I)QWC it may change with $(p+1)$-roots of
unity for the block of $\ker(A)$ in $A$ being $J_p$ even under
rigid motions
$(R,t)\in\mathbf{O}_{n+1}(\mathbb{C})\ltimes\mathbb{C}^{n+1}$
which preserve the canonical form, so it is unambiguous up to
$(p+1)$-roots of unity.

We have the diagonal Q(W)C respectively for
$A=\Sigma_{j=1}^{n+1}a_j^{-1}e_je_j^T,\
A=\Sigma_{j=1}^na_j^{-1}e_je_j^T$; the diagonal IQWC come in
different flavors, according to the block of $f_1:\
A=J_p+\Sigma_{j=p+1}^{n+1}a_j^{-1}e_je_j^T$; in particular if
$A=J_{n+1}$, then $\mathrm{spec}(A)=\{0\}$ is unambiguous. General
quadrics are those for which all eigenvalues have geometric
multiplicity $1$; equivalently each eigenvalue has an only
corresponding SJ block; in this case the quadric also admits
elliptic coordinates.

There are continuous groups of symmetries which preserve the SJ
canonical form for more than one SJ block corresponding to an
eigenvalue, so from a metric point of view a metric classification
according to the elliptic coordinates and continuous symmetries
may be a better one.

With $R_z:=I_{n+1}-zA,\
z\in\mathbb{C}\setminus\mathrm{spec}(A)^{-1}$ the family of
quadrics $\{x_z\}_z$ confocal to $x_0$ is given by
$Q_z(x_z)=x_z^TAR_z^{-1}x_z+2(R_z^{-1}B)^Tx_z+C+zB^TR_z^{-1}B=0$.
For $z\in\mathrm{spec}(A)^{-1}$ we obtain singular confocal
quadrics; those with $z^{-1}$ having geometric multiplicity $1$
admit a singular set which is an $(n-1)$-dimensional quadric
projectively equivalent to $C(\infty)$, so they will play an
important r\^{o}le in the discussion of homographies
$H\in\mathbf{PGL}_{n+1}(\mathbb{C})$ taking a confocal family into
another one, since $H^{-1}(C(\infty)),\ C(\infty)$ respectively
$C(\infty),\ H(C(\infty))$ will suffice to determine each confocal
family.

The Ivory affinity is an affine correspondence between confocal
quadrics and having good metric properties (it may be the reason
why Bianchi calls it {\it affinity} in more than one language): it
is given by $x_z=\sqrt{R_z}x_0+C(z),\
C(z):=-(\frac{1}{2}\int_0^z(\sqrt{R_w})^{-1}dw)B$. Note that
$C(z)=0$ for QC, $=\frac{z}{2}e_{n+1}$ for QWC; for IQWC it is the
Taylor series of $\frac{1}{2}\int_0^z(\sqrt{1-w})^{-1}dw$ at $z=0$
with each monomial $z^{k+1}$ replaced by $z^{k+1}J_p^k\bar f_1$,
where $J_p$ is the block of $f_1$ in $A$ and thus a polynomial of
degree $p$ in $z$. Note
$AC(z)+(I_{n+1}-\sqrt{R_z})B=0=(I_{n+1}+\sqrt{R_z})C(z)+zB$.
Applying $d$ to $Q_z(x_z)=0$ we get $dx_z^TR_z^{-1}(Ax_z+B)=0$, so
the unit normal $N_z$ is proportional to $\hat N_z:=-2\pa_zx_z$.
If $\mathbb{C}^{n+1}\ni x\in x_{z_1},x_{z_2}$, then $\hat
N_{z_j}=R_{z_j}^{-1}(Ax+B)$; using $R_z^{-1}-I_{n+1}=zAR_z^{-1},\
z_1R_{z_1}^{-1}-z_2R_{z_2}^{-1}=(z_1-z_2)R_{z_1}^{-1}R_{z_2}^{-1}$
we get $0=Q_{z_1}(x)-Q_{z_2}(x)=(z_1-z_2)\hat N_{z_1}^T\hat
N_{z_2}$, so two confocal quadrics cut each other orthogonally
(Lam\'{e}). For general quadrics the polynomial equation
$Q_z(x)=0$ has degree $n+1$ in $z$ and it has multiple roots iff
$0=\pa_zQ_z(x)=|\hat N_z|^2$; thus outside the locus of isotropic
normals elliptic coordinates (given by the roots $z_1,...,z_{n+1}$
of the said equation) give a parametrization of $\mathbb{C}^{n+1}$
suited to confocal quadrics.

We have now some classical metric properties of the Ivory
affinity: with $x_0^0,x_0^1\in x_0,\ V_0^1:=x_z^1-x_0^0$, etc the
Ivory Theorem (preservation of length of segments between confocal
quadrics) becomes
$|V_0^1|^2=|x_0^0+x_0^1-C(z)|^2-2(x_0^0)^T(I_{n+1}+\sqrt{R_z})x_0^1+zC=|V_1^0|^2$;
the preservation of lengths of rulings: $w_0^TAw_0=w_0^T\hat
N_0=0,\ w_z=\sqrt{R_z}w_0\Rightarrow
w_z^Tw_z=|w_0|^2-zw_0^TAw_0=|w_0|^2$; the symmetry of the {\it
tangency configuration} (TC): $(V_0^1)^T\hat N_0^0
=(x_0^0)^TA\sqrt{R_z}x_0^1-B^T(x_z^0+x_z^1-C(z))+C=(V_1^0)^T\hat
N_0^1$; the preservation of angles between segments and rulings:
$(V_0^1)^Tw_0^0+(V_1^0)^Tw_z^0=-z(\hat N_0^0)^Tw_0^0=0$; the
preservation of angles between rulings:
$(w_0^0)^Tw_z^1=(w_0^0)^T\sqrt{R_z}w_0^1=(w_z^0)^Tw_0^1$; the
preservation of angles between polar rulings: $(w_0^0)^TA\hat
w_0^0=0\Rightarrow (w_z^0)^T\hat w_z^0=(w_0^0)^T\hat
w_0^0-z(w_0^0)^TA\hat w_0^0=(w_0^0)^T\hat w_0^0$.

All complex quadrics are affine equivalent to either the unit
sphere $X\subset\mathbb{C}^{n+1},\ |X|^2=1$ or to the equilateral
paraboloid $Z\subset\mathbb{C}^{n+1},\ Z^T(I_{1,n}Z-2e_{n+1})=0$,
so a parametrization with regard to these two quadrics is in
order: $x_0=(\sqrt{A})^{-1}X$ for QC,
$x_0=(\sqrt{A+e_{n+1}e_{n+1}^T})^{-1}Z$ for QWC (for this reason
from a canonical metric point of view (that is we are interested
in a simplest form of $|dx_0|^2$) we should rather require that
$A^{-1}$ or $(A+e_{n+1}e_{n+1}^T)^{-1}$ is SJ).

For IQWC such a parametrization fails because of the isotropic
$\mathrm{ker}(A)$; however the computations between confocal
quadrics involving the Ivory affinity reveal a natural
parametrization of IQWC which is again an affine transformation of
$Z$.

Consider a canonical IQWC $x_0^T(Ax_0-2\bar f_1)=0,\
\ker(A)=\mathbb{C}f_1,\ A=J_p\oplus ...$ SJ. We are looking for a
linear map $L\in\mathbf{GL}_{n+1}(\mathbb{C})$ such that $x_0=LZ$,
equivalently $L^TAL=e^{2a}I_{1,n},\
I_{1,n}:=I_{n+1}-e_{n+1}e_{n+1}^T,\ L^T\bar f_1=e^{2a}e_{n+1}$.
Replacing $L$ with $L(e^{-a}I_{1,n}+e^{-2a}e_{n+1}e_{n+1}^T)$ we
can make $a=0$. Thus $Le_{n+1}=f_1,\ L^T(A+\bar f_1\bar
f_1^T)L=I_{n+1}$, so $L^{-1}=R^T\sqrt{A+\bar f_1\bar f_1^T},\
R^TR=I_{n+1}$ with $Re_{n+1}=\sqrt{A+\bar f_1\bar f_1^T}f_1$ (note
that $Re_{n+1}$ has, as required, length $1$). Once
$R\in\mathbf{O}_{n+1}(\mathbb{C})$ with the above property is
found, $L$ thus defined satisfies $L^T\bar f_1=e_{n+1}$ and thus
$L^TAL=I_{1,n}$. $L$ with the above properties is unique up to
rotations fixing $e_{n+1}$ in its domain and a canonical choice of
$R$ reveals itself from a SJ canonical form when doing
computations on confocal quadrics. We have $LL^T=(A+\bar f_1\bar
f_1^T)^{-1},\
I_{1,n}L^{-1}\sqrt{R_z}L=I_{1,n}L^{-1}\sqrt{R_z}LI_{1,n}=L^{-1}\sqrt{R_z}L-e_{n+1}\bar
f_1\sqrt{R_z}L=L^TA\sqrt{R_z}L=I_{1,n}\sqrt{I_{n+1}-zL^TA^2L}=:
I_{1,n}\sqrt{R'_z},\ A':=L^TA^2L,\
\ker(A')=\mathbb{C}e_{n+1}\oplus\mathbb{C}L^{-1}(A+\bar f_1\bar
f_1^T)^{-1}f_1 =\mathbb{C}e_{n+1}\oplus\mathbb{C}L^Tf_1$; choose
$R$ which makes $A'$ SJ. Note that we can take for QWC
$L:=(\sqrt{A+e_{n+1}e_{n+1}^T})^{-1},\ A':=A,\
\ker(A')=\mathbb{C}e_{n+1}$, so IQWC can be regarded as metrically
degenerated QWC. Note that
$e_{n+1}^TL^{-1}\sqrt{R_z}L=(-I_{1,n}L^{-1}C(z)+e_{n+1})^T$; this
can be confirmed analytically by differentiating with respect to
$z$ and using $(L^T)^{-1}=AL-Be_{n+1}^T$ and will imply the
symmetry of the TC, but since we have already proved the symmetry
of the TC, we can use this to imply the previous. Thus
$L^{-1}x_z=L^{-1}\sqrt{R_z}LZ+L^{-1}C(z)=I_{1,n}\sqrt{R'_z}Z+
e_{n+1}(-I_{1,n}L^{-1}C(z)+e_{n+1})^TZ+L^{-1}C(z),\
(x_z^1-x_0^0)^T\hat
N_0^0=(L^{-1}x_z^1-Z_0)^T(I_{1,n}Z_0-e_{n+1})=Z_0^TI_{1,n}\sqrt{R'_z}Z_1+
(Z_0+Z_1)^T(I_{1,n}L^{-1}C(z)-e_{n+1})-e_{n+1}^TL^{-1}C(z)$. Note
that for IQWC $|I_{1,n}L^{-1}C(z)|^2=2e_{n+1}^TL^{-1}C(z)$, so
$L^{-1}C(z)$ lies itself on $Z$ (also in this case since $\bar
f_1^TJ_p^k\bar f_1=\del_{k\ p-1}$ we have
$e_{n+1}^TL^{-1}C(z)=\bar
f_1^TC(z)=(^{-\frac{1}{2}}_{p-1})\frac{(-z)^p}{-2p}$, so
$e_{n+1}^T$ picks up the highest power of $z$ in $L^{-1}C(z)$). To
see this we need $0=|L^{-1}C(z)-\bar f_1^TC(z)e_{n+1}|^2-2\bar
f_1^TC(z)=C(z)^T(LL^T)^{-1}C(z)-(\bar f_1^TC(z))^2-2\bar
f_1^TC(z)=C(z)^TAC(z)-2\bar f_1^TC(z)$; using
$AC(z)=(I_{n+1}-\sqrt{R_z})\bar f_1$ and
$(I_{n+1}+\sqrt{R_z})C(z)=z\bar f_1$ it is satisfied.

Note also that as needed later we have
$(L^TL)^{-1}=A'-I_{1,n}L^{-1}Be_{n+1}^T-e_{n+1}(I_{1,n}L^{-1}B)^T+|B|^2e_{n+1}e_{n+1}^T,\
|\hat
N_0|^2=|(L^T)^{-1}(I_{1,n}Z-e_{n+1})|^2=Z^TA'Z+2Z^TI_{1,n}L^{-1}B+|B|^2,\\
(I_{n+1}+\sqrt{R'_z})I_{1,n}L^{-1}C(z)=I_{1,n}L^{-1}(I_{n+1}+\sqrt{R_z})C(z)=-zI_{1,n}L^{-1}B$.

\section{Multiply conjugate systems containing deformations of quadrics}

We assume all subspaces in discussion to be non-isotropic (the
Euclidean product induced on them by the one on $\mathbb{C}^m$ is
non-degenerate; this assures the existence of orthonormal normal
frames and thus the discussion of sub-manifolds via the {\it
Gau\ss-Weingarten} (GW) and {\it
Gau\ss-Codazzi-Mainardi(-Peterson)-Ricci} (G-CMP-R) equations).

For deformations $x\subset\mathbb{C}^{2n-1}$ of
$x_0\subset\mathbb{C}^{n+1}$ (that is $|dx|^2=|dx_0|^2$) with
common conjugate system $(u^1,...,u^n)$ and non-degenerate joined
second fundamental forms (that is $[d^2x_0^TN_0\ \ d^2x^TN]$ is a
symmetric quadratic $\mathbb{C}^n$-valued form which contains only
$(du^j)^2$ terms for $N_0$ unit normal field of $x_0$ and $N=[N_1\
\ ...\ \ N_{n-1}]$ orthonormal normal frame of $x$ and the
dimension $n$ cannot be lowered in an open dense set) the linear
element must satisfy the condition
$$\Ga_{jk}^l=0,\ j,k,l\ \mathrm{distinct}$$ and such deformations
$x\subset\mathbb{C}^{2n-1}$ are in bijective correspondence with
solutions $\{\mathbf{a}_j\}_{j=1,...,n}\subset\mathbf{C}^*$ of the
differential system
\begin{eqnarray}\label{eq:syst0}
(\log\mathbf{a}_j)_k=\Ga_{jk}^j,\ j\neq k,\
\sum_{j=1}^n\frac{(h_j^0)^2}{\mathbf{a}_j^2}+1=0,
\end{eqnarray}
where $N_0^Td^2x_0=:\sum_{j=1}^nh_j^0(du^j)^2$ is the second
fundamental form of $x_0$ (we shall use Latin indices $j,k,l,...$
including to differentiate respectively with $u^j,u^k,u^l,...$
when clear from the context; also we shall preserve the classical
notation $d^2$ for the symmetric (tensorial) second derivative and
we shall use $d\wedge$ for the exterior (antisymmetric)
derivative; thus $d\wedge d=0$).

Once a solution of this system is known, one finds the second
fundamental form of $x$ (complete the row
$[i\frac{h_1^0}{\mathbf{a}_1}\ \ ...\ \
i\frac{h_n^0}{\mathbf{a}_n}]$ as the first row in an orthogonal
$R\subset\mathbf{O}_n(\mathbb{C})$, delete it, multiply the column
$j$ of the $\mathbf{M}_{n-1,n}(\mathbb{C})$ obtained matrix
respectively with $\mathbf{a}_j$ and take the $k$-th row to obtain
the second fundamental form of $x$ in the $N_k$ direction) and
then one finds $x$ by the integration of a Ricatti equation and
quadratures (the Gau\ss-Bonnet(-Peterson) Theorem).

By the argument of Cartan's reduction of exteriorly orthogonal
forms to the canonical form such coordinates $\{u^j\}_j$ exist for
real deformations $\subset\mathbb{R}^{2n-1}$ of imaginary quadrics
$\subset\mathbb{R}^n\times i\mathbb{R}$ of negative curvature (for
such cases also all computations in the deformation problem will
be real; see Berger, Bryant, Griffiths \cite{BBG}). However, since
we have completely integrable differential systems (systems in
involution) for the deformation problem for quadrics, the
dimensionality of the space of deformations of quadrics remains
the same (namely solution depending on $n(n-1)$ functions of one
variable) in the complex setting also (the Cartan characters
remain the same).

To obtain multiply conjugate systems we need to extend $x$ with
the independent variables $u^{n+1},...,u^{2n-1}$.

An $m$-dimensional region
$$x=x(u^1,...,u^m)\subseteq\mathbb{C}^m,\ du^1\wedge...\wedge du^m\neq 0$$
gives a multiply conjugate system iff
\begin{eqnarray}\label{eq:multipl}
x_{jk}=(\log\mathbf{a}_j)_kx_j+(\log\mathbf{a}_k)_jx_k,\
\mathbf{a}_j\subset\mathbb{C}^*,\ j,k=1,...,m,\ j\neq k
\end{eqnarray}
with the compatibility condition
\begin{eqnarray}\label{eq:compm}
(\mathbf{a}_j)_{kl}=(\log\mathbf{a}_k)_l(\mathbf{a}_j)_k+(\log\mathbf{a}_l)_k(\mathbf{a}_j)_l,\
j,k,l=1,...,m\ \mathrm{distinct}.
\end{eqnarray}
Note that (\ref{eq:compm}) is the condition that the Riemann
symbols $R_{jklp}$ of the linear element
$\sum_{j=1}^m\mathbf{a}_j^2(du^j)^2$ with at least three of
$j,k,l,p$ distinct are $0$, so multiply conjugate systems in
$\mathbb{C}^m$ are put in correspondence with such linear elements
(normal curved spaces in Bianchi's denomination).

For the specific computations of deformations of quadrics we shall
use the convention $\mathbb{C}^n\subset\mathbb{C}^{n+1}$ with $0$
on the $(n+1)^{\mathrm{th}}$ component; thus for example we can
multiply $(n+1,n+1)$-matrices with $n$-column vectors and
similarly one can extend $(n,n)$ matrices to $(n+1,n+1)$ matrices
with zeroes on the last column and row. The converse is also
valid: an $(n+1,n+1)$ matrix with zeroes on the last column and
row (or multiplied on the left with an $n$-row vector and on the
right with an $n$-column vector) will be considered as an
$(n,n)$-matrix.

\subsection{(Isotropic) quadrics without center}\noindent

With $V:=\sum_{k=1}^nv^ke_k=[v^1\ \ ...\ \ v^n]^T$ consider the
complex equilateral paraboloid
$Z=Z(v^1,...,v^n)=V+\frac{|V|^2}{2}e_{n+1}$.

We have the (I)QWC $x_0:=LZ,\ L\in\mathbf{GL}_{n+1}(\mathbb{C})$
(recall $L:=(\sqrt{A+e_{n+1}e_{n+1}^T})^{-1},\
\ker(A)=\mathbb{C}e_{n+1}, A$ SJ, $B=-e_{n+1}$ for QWC and
$Le_{n+1}=f_1,\ L^T(A+\bar f_1\bar f_1^T)L=I_{n+1},\ A':=L^TA^2L$
SJ for $\ker(A)=\mathbb{C}f_1, A$ SJ, $B=-\bar f_1$ in the case of
IQWC) with linear element, unit normal, second fundamental form
and Christoffel symbols
$|dx_0|^2=dV^TL^TLdV+(V^TdV)^2|Le_{n+1}|^2+2(Le_{n+1})^TLdV(V^TdV),\
N_0=\frac{(L^T)^{-1}V+B}{\sqrt{H}},\
N_0^Td^2x_0=-\frac{|dV|^2}{\sqrt{H}},\
H:=|(L^T)^{-1}V+B|^2=V^TA'V+2V^TL^{-1}B+|B|^2,\ \ti\Ga_{jk}^l=0,\
j\neq k,\ \ti\Ga_{jj}^k=\frac{\pa\log\sqrt{H}}{\pa v^k}$.

Because $(v^1,...,v^n)$ are isothermal-conjugate and
$(u^1,..,u^n)$ are conjugate on $x_0$, the Jacobian
$\frac{\pa(v^1,...,v^n)}{\pa(u^1,...,u^n)}$ has orthogonal
columns, so with $\la_j:=|\frac{\pa V}{\pa u^j}|\neq 0,\
\La:=[\la_1\ \ ...\ \ \la_n]^T,\ \del':=\mathrm{diag}[du^1\ \ ...\
\ du^n],\ d'f:=\sum_{j=1}^nf_jdu^j$ we have $d'V=R\del'\La,\
R\subset\mathbf{O}_n(\mathbb{C})$. Multiplying the formula for the
change of Christoffel symbols $\frac{\pa v^c}{\pa
u^l}\Ga_{jk}^l=\frac{\pa^2v^c}{\pa u^j\pa u^k}+\frac{\pa v^a}{\pa
u^j}\frac{\pa v^b}{\pa u^k}\ti \Ga_{ab}^c=\frac{\pa^2v^c}{\pa
u^j\pa u^k}+\la_j^2\del_{jk}\frac{\pa\log\sqrt{H}}{\pa v^c}$ on
the left with $\frac{\pa v^c}{\pa u^p}$ and summing after $c$ we
obtain $\Ga_{jk}^p=\la_p^{-2}(\sum_c\frac{\pa v^c}{\pa
u^p}\frac{\pa^2v^c}{\pa u^j\pa
u^k}+\la_j^2\del_{jk}(\log\sqrt{H})_p)=\del_{pk}(\log\la_k)_j
+\del_{jk}\frac{\la_j^2}{\la_p^2}(\log\frac{\sqrt{H}}{\la_j})_p
+\del_{pj}(\log\la_j)_k$, so $\Ga_{jk}^j=(\log\la_j)_k,\
\Ga_{jj}^k=\frac{\la_j^2}{\la_k^2}(\log\frac{\sqrt{H}}{\la_j})_k,\
j\neq k,\ \Ga_{jj}^j=(\log(\la_j\sqrt{H}))_j$. We have
$h_j^0=-\frac{\la_j^2}{\sqrt{H}}$; since
$(\log\la_j)_k=\Ga_{jk}^j=(\log\mathbf{a}_j)_k,\ j\neq k$ we get
$\la_j=\phi_j(u^j)\mathbf{a}_j$; after a change of the $u^j$
variable into itself we can make $\la_j=\mathbf{a}_j,\ j=1,...,n$,
so from (\ref{eq:syst0}) $|\La|^2=-H$ and
$\La^Td'\La=-d'V^T(A'V+L^{-1}B)=-\La^T\del'R^T(A'V+L^{-1}B)$.

Imposing the compatibility condition $R^Td'\wedge$ on
$d'V=R\del'\La$ we get $R^Td'R\wedge\del'\La-\del'\wedge d'\La=0$,
or $(\la_j)_k=e_j^TR^TR_je_k\la_k,\ j\neq k,\
e_l^T\frac{R^TR_j}{\la_j}e_k=e_l^T\frac{R^TR_k}{\la_k}e_j,\ j,k,l$
distinct. Now by the standard Cartan trick
$-e_k^T\frac{R^TR_l}{\la_l}e_j=-e_k^T\frac{R^TR_j}{\la_j}e_l=e_l^T\frac{R^TR_j}{\la_j}e_k
=e_l^T\frac{R^TR_k}{\la_k}e_j=-e_j^T\frac{R^TR_k}{\la_k}e_l=-e_j^T\frac{R^TR_l}{\la_l}e_k
=e_k^T\frac{R^TR_l}{\la_l}e_j$ for $j,k,l$ distinct, so
$e_j^TR^TR_ke_l=0$ for $j,k,l$ distinct. Keeping account of the
prime integral property $\La^T[d'\La+\del'R^T(A'V+L^{-1}B)]=0$ and
with
$\om':=\sum_{j=1}^n(e_je_j^TR^TR_j\del'+\del'R^TR_je_je_j^T)=-\om'^T$
we have $d'\La=\om'\La-\del'R^T(A'V+L^{-1}B)$ and $d'\wedge
d'\La=0$ becomes the differential system
\begin{eqnarray}\label{eq:defqwc}
d'\wedge\om'=\om'\wedge\om'-\del'R^TA'R\wedge\del',\
\om'\wedge\del'=\del'\wedge R^Td'R,\
R\subset\mathbf{O}_n(\mathbb{C})
\end{eqnarray}
in involution (that is no further conditions appear if one imposes
$d'\wedge$ conditions and one uses the equations of the system
itself) as the compatibility condition for the completely
integrable linear differential system
\begin{eqnarray}\label{eq:systqwc}
d'V=R\del'\La,\
d'\La=\om'\La-\del'R^T(A'V+L^{-1}B),\ \La^T\La=-(V^TA'V+2V^TL^{-1}B+|B|^2).\nonumber\\
\end{eqnarray}
Note that for QWC once we know a solution
$R\subset\mathbf{O}_n(\mathbb{C})$ of (\ref{eq:defqwc}), a
solution $V,\ \La$ of (\ref{eq:systqwc}) and certain linearly
independent solutions $V_j,\ \La_j,\ j=1,...,2n-1$ of the
(homogeneous) differential part of (\ref{eq:systqwc}), then we can
find up to rigid motions the space realization of the deformation
$x\subset\mathbb{C}^{2n-1}$ of $x_0$ up to quadratures and without
the use of the Gau\ss-Bonnet(-Peterson) Theorem (this is again due
to Bianchi for $n=2$).

Let $V_j,\ \La_j,\ j=1,...,2n-1$ be certain linearly independent
solutions of the (homogeneous) differential part of
(\ref{eq:systqwc}) and $V,\ \La$ solution of (\ref{eq:systqwc})
such that $x\subset\mathbb{C}^{2n-1}$ given up to quadratures by
\begin{eqnarray}\label{eq:d'x}
d'x:=[V_1\ \ ...\ \ V_{2n-1}]^Td'V
\end{eqnarray}
is a deformation of $x_0$ with non-degenerate joined second
fundamental forms. Note $d'\wedge d'x=0$; since
$|d'x|^2=|d'x_0|^2,\ d'x_0=(L+e_{n+1}V^T)d'V$ we get
$\sum_{j=1}^{2n-1}V_jV_j^T-VV^T=I_{1,n}L^TL$; applying $d'$ we get
$\sum_{j=1}^{2n-1}V_j\La_j^T-V\La^T=0$ (use $\del'M=N\del',\
M,N\subset\mathbf{M}_n(\mathbb{C})\Leftrightarrow
M=N=\mathrm{diag}$); applying $d'$ again  we get
$\sum_{j=1}^{2n-1}\La_j\La_j^T-\La\La^T=I_n$.

With $\mathcal{V}:=[V_1\ \ ...\ \ V_{2n-1}\ \ iV],\
\mathcal{L}:=[\La_1\ \ ...\ \ \La_{2n-1}\ \ i\La]$ these can be
written as
\begin{eqnarray}\label{eq:vla0}
\begin{bmatrix}\mathcal{V}\\\mathcal{L}\end{bmatrix}[\mathcal{V}^T\
\ \mathcal{L}^T]=\begin{bmatrix}I_{1,n}L^TL& 0\\0&
I_n\end{bmatrix}.
\end{eqnarray}
We also have the prime integral property $[\mathcal{V}^T\ \
\mathcal{L}^T]\begin{bmatrix}A'&0\\0&I_n\end{bmatrix}
\begin{bmatrix}\mathcal{V}\\\mathcal{L}\end{bmatrix}=
\begin{bmatrix}\mathcal{C}&ic\\ic^T&1\end{bmatrix},\
\mathcal{C}=\mathcal{C}^T\in\mathbf{M}_{2n-1}(\mathbb{C}),\\
c\in\mathbb{C}^{2n-1}$; multiplying it on the left with
$\begin{bmatrix}\mathcal{V}\\\mathcal{L}\end{bmatrix}$ and using
(\ref{eq:vla0}) we get $\mathcal{C}=I_{2n-1},\ c=0$; thus
$[\mathcal{V}^T\ \ \mathcal{L}^T]$ is determined modulo a
multiplication on the left with a rotation
$\begin{bmatrix}\mathcal{R}&0\\0&1\end{bmatrix}\in\mathbf{O}_{2n}(\mathbb{C})$
by
\begin{eqnarray}\label{eq:vla}
[\mathcal{V}^T\ \
\mathcal{L}^T]\begin{bmatrix}A'&0\\0&I_n\end{bmatrix}
\begin{bmatrix}\mathcal{V}\\\mathcal{L}\end{bmatrix}=I_{2n}.
\end{eqnarray}
Thus  (\ref{eq:vla}) is equivalent to (\ref{eq:vla0}) and to
$[\mathcal{V}^TL^{-1}\ \
\mathcal{L}^T]\subset\mathbf{O}_{2n}(\mathbb{C})$.

For the non-degenerate joined second fundamental forms property we
need to prove that there is no vector field
$N\subset\mathbb{C}^{2n-1}$ along $x$ such that $N^Td'x=0,\
N^Td'^2x=|d'V|^2=\La^T\del'^2\La$, that is the linear system
$\begin{bmatrix}\mathcal{V}\\\mathcal{L}\end{bmatrix}\begin{bmatrix}N\\0\end{bmatrix}
=\begin{bmatrix}0\\\La\end{bmatrix}$ is inconsistent; using
(\ref{eq:vla}) this becomes $\begin{bmatrix}N\\0\end{bmatrix}
=\mathcal{L}^T\La$, which is indeed inconsistent because
$|\La|^2\neq 0$.

Note that since for IQWC
$\ker{A'}=\mathbb{C}e_{n+1}\oplus\mathbb{C}L^Tf_1,\
\mathrm{coker}{A'}=\mathbb{C}e_{n+1}\oplus\mathbb{C}L^{-1}B$ for
such quadrics we cannot derive from the equivalent of
(\ref{eq:vla}) for IQWC the full information about the equivalent
of (\ref{eq:vla0}) for IQWC. This is to be expected, since a
fundamental set of solutions of the homogeneous differential part
of (\ref{eq:systqwc}) contains $2n$ linearly independent
solutions; the obvious solution $[f_1^TL\ \ 0]^T$ is isotropic and
perpendicular on all the others {\it with respect to} (wrt) the
bilinear form given by the left hand side of (\ref{eq:vla}); by
the Gramm-Schmidt orthogonalization process one can complete this
solution with $2n-1$ other solutions $V_j,\ \La_j,\ j=1,...,2n-1$
orthonormal wrt the same bilinear form. This isotropic solution is
added in the $2n-1$ other with yet undetermined constant
coefficients $V_j^TL^{-1}B$. The part of the equivalent of
(\ref{eq:vla0}) for IQWC that cannot be inferred from the
equivalent of (\ref{eq:vla}) for IQWC is
$\begin{bmatrix}V\\\La\end{bmatrix}=-\sum_{j=1}^{2n-1}
\begin{bmatrix}V_j\\\La_j\end{bmatrix}V_j^TL^{-1}B+
\begin{bmatrix}L^Tf_1\\0\end{bmatrix}V^TL^{-1}B$, but this
requires $V^TL^{-1}B$ to be constant, so $V,\ \La$ has to be
solution of the homogeneous differential part of
(\ref{eq:systqwc}), a contradiction, so $V_j,\ \La_j,\
j=1,...,2n-1$ cannot be solutions of the homogeneous differential
part of (\ref{eq:systqwc}).

We now extend the independent variables $u^1,...,u^n$ with
independent variables $u^{n+1},...,u^{2n-1},\\ du^1\wedge...\wedge
du^{2n-1}\neq 0,\ d=d'+d'',\ d'f=\sum_{j=1}^nf_jdu^j,\
d''f=\sum_{j=1}^{n-1}f_{j+n}du^{j+n}$; thus $d\wedge d=0$ becomes
$d'\wedge d'=0,\ d'\wedge d''+d''\wedge d'=0,\ d''\wedge d''=0$.

With $\del'':=\mathrm{diag}[du^{n+1}\ \ ...\ \ du^{2n-1}\ \ 0]$ we
now have $d''\La=M\del''V,\ M\subset\mathbf{M}_n(\mathbb{C}),\
Me_n=0$ (this choice of $M$ is due to Bianchi for $n=2$ and its
reason will appear immediately) .

Imposing the compatibility condition $d'\wedge d''\La+d''\wedge
d'\La=0$ we obtain $[d'\wedge(M\del'')-\om'\wedge
M\del'']V+(d''\wedge\om'-M\del''\wedge R\del')\La
+\del'R^T\wedge(A'd''V-d''RR^T(A'V+L^{-1}B))=0$, from where we
obtain by applying $\del'\wedge:\
\del'\wedge[d'\wedge(M\del'')-\om'\wedge M\del'']=0,\
\del'\wedge(d''\wedge\om'-M\del''\wedge R\del')=0$, so
$d'\wedge(M\del'')=\om'\wedge M\del''+\del'R^T\wedge N\del'',\
N\subset\mathbf{M}_n(\mathbb{C}),\ Ne_n=0,\
d''\wedge\om'=M\del''\wedge R\del'+\del'R^T\wedge ...$. Since
$\om'^T=-\om'$ we have $d''\wedge\om'=M\del''\wedge
R\del'+\del'R^T\wedge\del''M^T$ (in particular we get $M$ from $R$
and its derivatives; this is the reason for the choice of $M$), so
$A'd''V=(d''RR^TA'-N\del'')V-\del''M^T\La+d''RR^TL^{-1}B.$
Imposing the compatibility condition $d'\wedge A'd''V+d''\wedge
A'd'V=0$ and collecting the coefficient of $R\del'\La$ we get
$d''RR^TA'-N\del''=A'd''RR^T+\del''N^T$ (in particular we get $N$
from $R$ and its derivatives). Note that since for IQWC
$\ker{A'}=\mathbb{C}e_{n+1}\oplus\mathbb{C}L^Tf_1,\
\mathrm{coker}{A'}=\mathbb{C}e_{n+1}\oplus\mathbb{C}L^{-1}B$ for
such quadrics we cannot derive the full information about $d''V$
from $A'd''V$ only: from the prime integral property
$d''\La^T\La=-d''V^T(A'V+L^{-1}B)$ and using
$V^T(d''RR^TA'-N\del'')V=0$ we have
$(L^{-1}B)^T(d''V-d''RR^TV)=0$. Now using
$A'(d''V-d''RR^TV)=[\mathcal{A}^{-1}+e_{n+1}(L^{-1}B)^T](d''V-d''RR^TV),\
\mathcal{A}:=L^TL$ we can finally extend the completely integrable
linear differential system (\ref{eq:systqwc}) to the completely
integrable linear differential system
\begin{eqnarray}\label{eq:systqwcext}
dV=R\del'\La+(d''RR^T+\mathcal{A}\del''N^T)V-\mathcal{A}\del''M^T\La
+\mathcal{A}d''RR^TL^{-1}B,\nonumber\\
d\La=\om'\La-\del'R^T(A'V+L^{-1}B)+M\del''V,\
\La^T\La=-(V^TA'V+2V^TL^{-1}B+|B|^2)
\end{eqnarray}
with the extended compatibility (and algebraic) conditions
\begin{eqnarray}\label{eq:defqwcext}
d'\wedge\om'=\om'\wedge\om'-\del'R^TA'R\wedge\del',\
\om'\wedge\del'=\del'\wedge R^Td'R,\
d''\wedge\om'=M\del''\wedge R\del'+\del'R^T\wedge\del''M^T\nonumber\\
(\Leftrightarrow^{Rd''\wedge(R^Td'R)=-d'\wedge(d''RR^T)R}d'\wedge(d''RR^T)
=\del''M^T\wedge\del'R^T+R\del'\wedge M\del''),\nonumber\\
d'\wedge(M\del'')=\om'\wedge M\del''+\del'R^T\wedge N\del'',\
d'\wedge(N\del'')=R\del'\wedge M\del''A'-A'R\del'\wedge
M\del'',\nonumber\\
d''\wedge(M\del'')=M\del''\wedge(d''RR^T+\mathcal{A}\del''N^T),\
d''\wedge(N\del'')=\del''M^T\wedge M\del''+\nonumber\\
N\del''\wedge(d''RR^T+\mathcal{A}\del''N^T)+d''RR^T\wedge
N\del'',\ A'd''RR^T+\del''N^T=d''RR^TA'-N\del'',\nonumber\\
M\del''\mathcal{A}\wedge\del''M^T=0,\
e_{n+1}^T\mathcal{A}\del''M^T=e_{n+1}^T\mathcal{A}\del''N^T=0,\
e_{n+1}^T\mathcal{A}d''RR^TL^{-1}B=0,\nonumber\\
M\del''L^{-1}B=M\del''\mathcal{A}\wedge
d''RR^TL^{-1}B=N\del''\mathcal{A}\wedge d''RR^TL^{-1}B=0 ,\
R\subset\mathbf{O}_n(\mathbb{C}).
\end{eqnarray}
If we impose $d'\wedge$ and $d''\wedge$ on (\ref{eq:defqwcext}),
use $d'\wedge d'=0,\ d'\wedge d''+d''\wedge d'=0,\ d''\wedge
d''=0$ and the equations of the system itself, then we get the
algebraic conditions
\begin{eqnarray}\label{eq:newcond}
M\del''\mathcal{A}\wedge\del''N^T\wedge R\del'+\del'R^T\wedge
N\del''\wedge\mathcal{A}\del''M^T=0(\Leftrightarrow
M\del''\mathcal{A}\wedge\del''N^TR=\mathrm{diag}),\nonumber\\
M\del''A'^kL^{-1}B=N\del''A'^kL^{-1}B=M\del''\mathcal{A}\wedge
d''RR^TA'^kL^{-1}B=N\del''\mathcal{A}\wedge
d''RR^TA'^kL^{-1}B=0,\nonumber\\
k=0,1,...,p-1,\ A'=...\oplus J_p\oplus...\oplus J_1,\ p\ge 0\
\mathrm{being\ the \ SJ\ decomposition\ of}\ A'
\end{eqnarray}
(the last relations are relevant for IQWC and state that
$M\del'',\ N\del'',\ M\del''\mathcal{A}\wedge d''RR^T,\
N\del''\mathcal{A}\wedge d''RR^T$ are not supported on the
coordinates corresponding to the SJ block of $L^Tf_1$ in $A'$).

Further imposing $d'\wedge$ and $d''\wedge$ on (\ref{eq:newcond}),
using $d'\wedge d'=0,\ d'\wedge d''+d''\wedge d'=0,\ d''\wedge
d''=0$, the equations of (\ref{eq:newcond}) itself and
(\ref{eq:defqwcext}) we finally get involution (that is no further
conditions appear); note that for diagonal QWC (which form an open
dense set in the set of (I)QWC) (\ref{eq:newcond}) is vacuous, so
we have (\ref{eq:defqwcext}) already in involution.

While the differential system (\ref{eq:systqwcext}) together with
its compatibility conditions
(\ref{eq:defqwcext})\&(\ref{eq:newcond}) is interesting on its
own, these compatibility conditions must further be extended with
new conditions in order to describe multiply conjugate systems
containing deformations of quadrics; in order to find these
conditions we need to consider the space realization of solutions.

To extend $x$  with the independent variables
$u^{n+1},...,u^{2n-1}$ (\ref{eq:multipl}), (\ref{eq:compm}) and
(\ref{eq:d'x}) will provide the needed information to obtain
multiply conjugate systems containing deformations of QWC. For our
problem we have $m=2n-1,\ d'x=[V_1\ \ ...\ \ V_{2n-1}]^Td'V$ with
$R$ solution of (\ref{eq:defqwcext})\&(\ref{eq:newcond}), $V,\
\La,\ V_j,\ \La_j,\ j=1,...,2n-1$ solutions of the (homogeneous)
differential part of (\ref{eq:systqwcext}) satisfying
(\ref{eq:vla}) and $\mathbf{a}_j=\la_j,\ j=1,...,n$; from $x_{j\
n+l}=(\log\la_j)_{n+l}x_j+(\log\mathbf{a}_{n+l})_jx_{n+l},\
j=1,...,n,\ l=1,...,n-1$ we get $x_{n+l}=([V_1\ \
...V_{2n-1}]^TN-[\La_1\ \ ...\ \ \La_{2n-1}]^TM)e_l\frac{1}{c_l},\
d'(\log\mathbf{a}_{n+l})=c_le_l^T\mathcal{A}d'V$; from
$(\mathbf{a}_{n+l})_{jk}=(\log\la_j)_k(\mathbf{a}_{n+l})_j+(\log\la_k)_j(\mathbf{a}_{n+l})_k,\
j,k=1,...,n,\ j\neq k,\ l=1,...,n-1$ we get
$\frac{1}{c_l}=e_l^T(\mathcal{A}V+f),\ d'f=0$, so
$\mathbf{a}_{n+l}=e^{g_l}e_l^T(\mathcal{A}V+f),\ d'g_l=0,\
d''x=([V_1\ \ ...V_{2n-1}]^TN\del''-[\La_1\ \ ...\ \
\La_{2n-1}]^TM\del'')(\mathcal{A}V+f)$. We need $0=d'\wedge
d''x+d''\wedge d'x=[V_1\ \ ...\ \ V_{2n-1}]^TR\del'\wedge
M\del''A'f$ (so $f=0$), $0=d''\wedge d''x=([V_1\ \ ...\ \
V_{2n-1}]^TN\del''-[\La_1\ \ ...\ \
\La_{2n-1}]^TM\del'')\mathcal{A}^2\wedge(\del''M^T\La-
\del''N^TV)$, so we need
\begin{eqnarray}\label{eq:geomcond}
M\del''\mathcal{A}\wedge\mathcal{A}\del''M^T=
N\del''\mathcal{A}\wedge\mathcal{A}\del''M^T=
N\del''\mathcal{A}\wedge\mathcal{A}\del''N^T=0.
\end{eqnarray}
Imposing $d'\wedge$ and $d''\wedge$ conditions on
(\ref{eq:geomcond}) and using
(\ref{eq:defqwcext})\&(\ref{eq:newcond}) and (\ref{eq:geomcond})
themselves we don't get any new conditions; note also that
(\ref{eq:geomcond}) are vacuous for diagonal QWC. Now we need
\begin{eqnarray}\label{eq:cond}
x_{n+l\ n+p}=(\log\mathbf{a}_{n+l})_{n+p}x_{n+l}
+(\log\mathbf{a}_{n+p})_{n+l}x_{n+p},\ l,p=1,...,n-1,\ l\neq
p,\nonumber\\ (\la_j)_{k\
n+l}=(\log\la_k)_{n+l}(\la_j)_k+(\log\mathbf{a}_{n+l})_k(\la_j)_{n+l},\
j,k=1,...,n,\
j\neq k,\ l=1,...,n-1,\nonumber\\
(\la_j)_{n+l\ n+p}=(\log\mathbf{a}_{n+l})_{n+p}(\la_j)_{n+l}
+(\log\mathbf{a}_{n+p})_{n+l}(\la_j)_{n+p},\ j=1,...,n,\
l,p=1,...,n-1,\ l\neq p,\nonumber\\
(\mathbf{a}_{n+l})_{j\
n+p}=(\log\la_j)_{n+p}(\mathbf{a}_{n+l})_j+(\log\mathbf{a}_{n+p})_j(\mathbf{a}_{n+l})_{n+p},\
j=1,...,n,\
l,p=1,...,n-1,\ l\neq p,\nonumber\\
(\mathbf{a}_{n+l})_{n+p\
n+q}=(\log\mathbf{a}_{n+p})_{n+q}(\mathbf{a}_{n+l})_{n+p}
+(\log\mathbf{a}_{n+q})_{n+p}(\mathbf{a}_{n+l})_{n+q},\
l,p,q=1,...,n-1\ \mathrm{distinct};\nonumber\\
\end{eqnarray}
these are satisfied for diagonal QWC if
$d''g_l=(g_l)_{n+l}du^{n+l}$.

\subsection{Quadrics with center}\noindent

\section{The B\"{a}cklund transformation}

\subsection{(Isotropic) quadrics without center}\noindent

For the B transformation the space realization
\begin{eqnarray}\label{eq:x1qwc}
x^1=x^0+[x_{v_0^1}^0\ \ ...\ \
x_{v_0^n}^0](\sqrt{R'_z}V_1-V_0+I_{1,n}L^{-1}C(z))\subset\mathbb{C}^{2n-1}
\end{eqnarray}
of the leaf $x^1$ relative to the seed
$x^0\subset\mathbb{C}^{2n-1},\ |d'x^0|^2=|d'x_0^0|^2$, the
algebraic transformation
\begin{eqnarray}\label{eq:algqwc}
V_1=\sqrt{R'_z}V_0-\sqrt{z}R_1\La_0+I_{1,n}L^{-1}C(z),\
\La_1=R_0^T(\sqrt{z}A'V_0+\sqrt{R'_z}R_1\La_0+\sqrt{z}I_{1,n}L^{-1}B),\nonumber\\
(0,\sqrt{z})\leftrightarrow(1,-\sqrt{z})
\end{eqnarray}
of solutions of (\ref{eq:systqwcext}) and the algebraic formula of
the BPT
\begin{eqnarray}\label{eq:bptqwc}
R_3R_0^T=(D_2-D_1R_2R_1^T) (D_2R_2R_1^T-D_1)^{-1},\
D_j:=\sqrt{R'_{z_j}}/\sqrt{z_j},\ j=1,2
\end{eqnarray}
remain valid, but the differential system subjacent to the B
transformation (Ricatti equation)
\begin{eqnarray}\label{eq:bqwc}
-d'R_1=R_1\om'_0+R_1\del'R_0^TDR_1-DR_0\del',\
 D:=\sqrt{R'_z}/\sqrt{z},\ (0,\sqrt{z})\leftrightarrow(1,-\sqrt{z})
\end{eqnarray}
in $R_1$ must be extended and then the BPT algebraic formula must
satisfy this extension (for the third M\"{o}bius configuration the
algebraic computations only suffice).

Applying $d''$ to (\ref{eq:algqwc}) we get the full Ricatti
equation
\begin{eqnarray}\label{eq:d''bqwc}
d'R_1=-R_1\om'_0-R_1\del'R_0^TDR_1+DR_0\del',\nonumber\\
d''R_1=D^{-1}(d''R_0R_0^TD^2+N_0\del'')D^{-1}R_1+
D^{-1}\del''M_0^T-R_1M_0\del''D^{-1}R_1
\end{eqnarray}
in $R_1$ and further
\begin{eqnarray}\label{eq:d''algqwc}
M_1\del''=R_0^T(DR_1M_0\del''-N_0\del'')D^{-1},\
N_1\del''=(A'R_1M_0\del''-R_1M_0\del''A'+DN_0\del'')D^{-1}\nonumber\\
\end{eqnarray}
(thus we need $A'e_n\subset\mathbb{C}e_n$).

Imposing compatibility conditions $d'\wedge,\ d''\wedge$ on
(\ref{eq:d''bqwc}) and using the equation itself we get $d'\wedge
d'R_1=R_1(-d'\wedge\om'_0+\om'_0\wedge\om'_0-\del'R_0^TA'R_0\wedge\del')
-DR_0(\del'\wedge\om'_0-R_0^Td'R_0\wedge\del')
+R_1(\om_0'\wedge\del'-\del'\wedge R_0^Td'R_0)R_0^TDR_1,\ d'\wedge
d''R_1+d''\wedge
d'R_1=D^{-1}[d'\wedge(d''R_0R_0^T)-\del''M_0^T\wedge\del'R_0^T-R_0\del'\wedge
M_0\del'']DR_1-R_1[d''\wedge\om'_0-M_0\del''\wedge
R_0\del'-\del'R_0^T\wedge\del''M_0^T]-R_1[d'\wedge(M_0\del'')-\om'_0\wedge
M_0\del''-\del'R_0^T\wedge N_0\del'']D^{-1}R_1
+D^{-1}[d'\wedge(N_0\del'')-R_0\del'\wedge
M_0\del''A'+A'R_0\del'\wedge M_0\del'']D^{-1}R_1
+D^{-1}[d'\wedge(\del''M_0^T)-\del''M_0^T\wedge\om'_0+\del''N_0^T\wedge
R_0\del'],\ d''\wedge
d''R_1=D^{-1}[d''\wedge(N_0\del'')-\del''M_0^T\wedge
M_0\del''-N_0\del''\wedge(d''R_0R_0^T-D^{-2}\del''N_0^T)-d''R_0R_0^T\wedge
N_0\del'']D^{-1}R_1+D^{-1}[d''\wedge(\del''M_0^T)-(d''R_0R_0^T+N_0\del''D^{-2})\wedge\del''M_0^T]
-R_1[d''\wedge(M_0\del'')-M_0\del''\wedge(d''R_0R_0^T-D^{-2}\del''N_0^T)]
+R_1M_0\del''D^{-2}\wedge\del''M_0^T+(R_1M_0\del''-D^{-1}N_0\del'')D^{-2}\wedge
(d''R_0R_0^TD^2+N_0\del''-D^2d''R_0R_0^T+\del''N_0^T)D^{-1}R_1$;
for diagonal QWC these are $0$ from (\ref{eq:defqwcext}), so
(\ref{eq:d''bqwc}) is completely integrable and admits solution
for any initial value $R_1$. If the initial value is orthogonal,
we would like the solution to remain orthogonal: $d(R_1R_1^T-I_n)
=-R_1\del'R_0^TD(R_1R_1^T-I_n)-(R_1R_1^T-I_n)DR_0\del'R_1^T
-(R_1R_1^T-I_n)D^{-1}(\del''M_0^TR_1^T+d''R_0R_0^TD^2+N_0\del'')D^{-1}$,
so $R_1R_1^T-I_n$ is a solution of a linear differential equation
and remains $0$ if initially it was $0$. The fact that $R_1$ is
itself a solution of (\ref{eq:defqwcext}) follows from the
symmetry $(0,\sqrt{z})\leftrightarrow(1,-\sqrt{z})$ and the fact
that $d\wedge dR_0=0$ (basically we use the converse of the proven
results).

Therefore we only need to prove that $R_3$ given by
(\ref{eq:bptqwc}) satisfies (\ref{eq:d''bqwc}) for $(R_0,z)$
replaced by $(R_1,z_2),\ (R_2,z_1)$; by symmetry it is enough to
prove only one relation. Since
$dR_1=-R_1\om'_0-R_1\del'R_0^TD_1R_1+D_1R_0\del'
+D_1^{-1}(d''R_0R_0^TD_1^2+N_0\del'')D_1^{-1}R_1+
D_1^{-1}\del''M_0^T-R_1M_0\del''D_1^{-1}R_1,\
dR_2=-R_2\om'_0-R_2\del'R_0^TD_2R_2+D_2R_0\del'+D_2^{-1}(d''R_0R_0^TD_2^2+N_0\del'')D_2^{-1}R_2+
D_2^{-1}\del''M_0^T-R_2M_0\del''D_2^{-1}R_2$, we get
$d(R_2R_1^T)=[-R_2\om'_0-R_2\del'R_0^TD_2R_2+D_2R_0\del'
+D_2^{-1}(d''R_0R_0^TD_2^2+N_0\del'')D_2^{-1}R_2+
D_2^{-1}\del''M_0^T-R_2M_0\del''D_2^{-1}R_2]R_1^T
+R_2[\om'_0R_1^T-R_1^TD_1R_0\del'R_1^T+\del'R_0^TD_1
+R_1^TD_1^{-1}(-D_1^2d''R_0R_0^T+\del''N_0^T)D_1^{-1}+
M_0\del''D_1^{-1}-R_1^TD_1^{-1}\del''M_0^TR_1^T]=-(R_2R_1^T)R_1\del'R_0^T(D_2R_2R_1^T-D_1)
+(D_2-R_2R_1^TD_1)R_0\del'R_1^T
-(R_2R_1^TD_1^{-1}-D_2^{-1})\del''M_0^TR_1^T-R_2M_0\del''(D_2^{-1}R_2R_1^T-D_1^{-1})
+D_2^{-1}(d''R_0R_0^TD_2^2+N_0\del'')D_2^{-1}R_2R_1^T
-R_2R_1^TD_1^{-1}(d''R_0R_0^TD_1^2+N_0\del'')D_1^{-1}$. Thus if we
prove the similar relation
$d(R_3R_0^T)=\\-(R_3R_0^T)R_0\del'R_1^T(D_2R_3R_0^T+D_1)
+(D_2+R_3R_0^TD_1)R_1\del'R_0^T+(R_3R_0^TD_1^{-1}+D_2^{-1})\del''M_1^TR_0^T
-R_3M_1\del''(D_2^{-1}R_3R_0^T+D_1^{-1})+D_2^{-1}(d''R_1R_1^TD_2^2+N_1\del'')D_2^{-1}R_3R_0^T
-R_3R_0^TD_1^{-1}(d''R_1R_1^TD_1^2+N_1\del'')D_1^{-1}$, then since
$dR_0=-R_0\om_1+R_0\del'R_1^TD_1R_0-D_1R_1\del'+D_1^{-1}(d''R_1R_1^TD_1^2+N_1\del'')D_1^{-1}R_0-
D_1^{-1}\del''M_1^T+R_0M_1\del''D_1^{-1}R_0$ we obtain what we
want:
$dR_3=-R_3\om_1-R_3\del'R_1^TD_2R_3+D_2R_1\del'+D_2^{-1}(d''R_1R_1^TD_2^2+N_1\del'')D_2^{-1}R_3+
D_2^{-1}\del''M_1^T-R_3M_1\del''D_2^{-1}R_3$. Differentiating
(\ref{eq:bptqwc}) we get $d(R_3R_0^T)(D_2R_2R_1^T-D_1)=
-(R_3R_0^TD_2+D_1)d(R_2R_1^T)$; thus we need to prove
$[-(R_3R_0^T)R_0\del'R_1^T(D_2R_3R_0^T+D_1)
+(D_2+R_3R_0^TD_1)R_1\del'R_0^T+(R_3R_0^TD_1^{-1}+D_2^{-1})\del''M_1^TR_0^T
-R_3M_1\del''(D_2^{-1}R_3R_0^T+D_1^{-1})+D_2^{-1}(d''R_1R_1^TD_2^2+N_1\del'')D_2^{-1}R_3R_0^T
-R_3R_0^TD_1^{-1}(d''R_1R_1^TD_1^2+N_1\del'')D_1^{-1}](D_2R_2R_1^T-D_1)=
-(R_3R_0^TD_2+D_1)\\\ [-(R_2R_1^T)R_1\del'R_0^T(D_2R_2R_1^T-D_1)
+(D_2-R_2R_1^TD_1)R_0\del'R_1^T-(R_2R_1^TD_1^{-1}-D_2^{-1})\del''M_0^TR_1^T
-R_2M_0\del''(D_2^{-1}R_2R_1^T-D_1^{-1})+D_2^{-1}(d''R_0R_0^TD_2^2+N_0\del'')D_2^{-1}R_2R_1^T
-R_2R_1^TD_1^{-1}(d''R_0R_0^TD_1^2+N_0\del'')D_1^{-1}]$. The terms
containing $R_1\del'R_0^T$ become $D_2+R_3R_0^TD_1=
(R_3R_0^TD_2+D_1)R_2R_1^T$ which follows directly from
(\ref{eq:bptqwc}) and the terms containing $R_0\del'R_1^T$ become
$R_0\del R_1^T(D_2R_3R_0^T+D_1) (D_2R_2R_1^T-D_1)=(R_3R_0^T)^T
(R_3R_0^TD_2+D_1) (D_2-R_2R_1^TD_1)R_0\del R_1^T$ which follows
from $(D_2R_3R_0^T+D_1) (D_2R_2R_1^T-D_1)=D_2^2-D_1^2=
(\frac{1}{z_2}-\frac{1}{z_1})I_n=(R_3R_0^T)^T(R_3R_0^TD_2+D_1)
(D_2-R_2R_1^TD_1)$.

Replacing $d''R_1R_1^T,\ M_1\del'',\ N_1\del''$ from
(\ref{eq:d''bqwc})\&(\ref{eq:d''algqwc}) and using
$\del''N_0^T=-A'd''R_0R_0^T+d''R_0R_0^TA'-N\del''$ the remaining
terms split into ones containing $d''R_0R_0^T,\ R_1M_0\del'',\
\del''M_0^TR_1^T,\ N_0\del''$ and they turn out to be as they
should.

For the terms containing $\Om:=d''R_0R_0^T$ we need
$[D_2^{-1}D_1^{-1}\Om D_1D_2R_3R_0^T-R_3R_0^TD_1^{-2}\Om D_1^2
+(R_3R_0^TD_1^{-1}+D_2^{-1})D_1^{-1}(A'\Om-\Om
A')](D_2R_2R_1^T-D_1)+(R_3R_0^TD_2+D_1) (D_2^{-1}\Om
D_2R_2R_1^T\\-R_2R_1^TD_1^{-1}\Om D_1)=0$; for the terms
containing $\Om:=R_1M_0\del''$ we need $[D_2^{-1}(-\Om
D_2+(A'\Om-\Om
A')D_2^{-1})D_1^{-1}R_3R_0^T+R_3R_0^TD_1^{-1}(\Om-(A'\Om-\Om
A')D_1^{-2}) -R_3R_0^TD_1\Om
D_1^{-1}(D_2^{-1}R_3R_0^T-D_1^{-1})](D_2R_2R_1^T-D_1)-(R_3R_0^TD_2+D_1)
R_2R_1^T\Om(D_2^{-1}R_2R_1^T-D_1^{-1})=0$; for the terms
containing $\Om:=D_2^{-1}D_1^{-1}\del''M_0^TR_1^T$ we need
$\Om(D_2R_3R_0^T+D_1)(D_2R_2R_1^T-D_1)-(R_3R_0^TD_2+D_1)(R_2R_1^TD_2-D_1)\Om=0$;
for the terms containing $\Om:=N_0\del''$ we need
$[D_2^{-1}D_1^{-1}(\Om D_2^2+D_1^2\Om
)D_2^{-1}D_1^{-1}R_3R_0^T-R_3R_0^T(D_1^{-2}\Om+\Om
D_1^{-2})+(R_3R_0^TD_1^{-1}+D_2^{-1})D_1^{-1}\Om +R_3R_0^T\Om
D_1^{-1}(D_2^{-1}R_3R_0^T-D_1^{-1})](D_2R_2R_1^T-D_1)+(R_3R_0^TD_2+D_1)
[D_2^{-1}\Om D_2^{-1}R_2R_1^T-R_2R_1^TD_1^{-1}\Om D_1^{-1}]=0$.

\subsection{Quadrics with center}

\end{document}